# ON LINEAR AND UNCONDITIONALLY ENERGY STABLE ALGORITHMS FOR VARIABLE MOBILITY CAHN-HILLIARD TYPE EQUATION WITH LOGARITHMIC FLORY-HUGGINS POTENTIAL

XIAOFENG YANG † AND JIA ZHAO ‡

ABSTRACT. In this paper, we consider the numerical approximations for the fourth order Cahn-Hilliard equation with concentration dependent mobility, and the logarithmic Flory-Huggins potential. One challenge in solving such a diffusive system numerically is how to develop proper temporal discretization for nonlinear terms in order to preserve the energy stability at the time-discrete level. We resolve this issue by developing a set of the first and second order time marching schemes based on a novel, called "Invariant Energy Quadratization" approach. Its novelty is that the proposed scheme is linear and symmetric positive definite because all nonlinear terms are treated semi-explicitly. We further prove all proposed schemes are unconditionally energy stable rigorously. Various 2D and 3D numerical simulations are presented to demonstrate the stability, accuracy and efficiency of the proposed schemes thereafter.

## 1. INTRODUCTION

In this paper, we consider the numerical approximations for the fourth order Cahn-Hilliard equation with the variable concentration-dependent mobility and the logarithmic Flory-Huggins (F-H) bulk potential. The Cahn-Hilliard equation is one of the two typical equations of the diffusive interface approach, or referred as the phase field method. With a long trackable history dated back to Rayleigh [30] and Van der Waals [41] a century ago, this method has now become a well-known efficient modeling and numerical tool to resolve the motion of free interfaces between multiple material components. About the most recent advances in modeling, algorithms or computational technologies, we refer to [4, 6, 7, 10, 15, 20–23, 26, 28, 40, 49] and the references therein.

There exist many advantages in the phase field model from the mathematical point of view. In particular, the model is usually energy stable (thermodynamics-consistent) and well-posed due to the energy-based variational approach, making it possible to carry out effective mathematical/numerical analysis and perform reliable computer simulations. For algorithms design, a significant goal is to verify the energy stable property at the discrete level irrespectively of the coarseness of the discretization (in what follows, those algorithms will be called unconditionally energy stable or thermodynamically consistent). Schemes with this property is specially preferred for solving diffusive systems since it is not only critical for the numerical scheme to capture the correct long time dynamics of the system, but also provides sufficient flexibility for dealing with the stiffness issue. Meanwhile, since the dynamics of coarse-graining (macroscopic) process may undergo rapid changes near the interface, the noncompliance of energy dissipation laws may lead to spurious numerical solutions if the mesh or time step size are not carefully controlled. We also emphasize that the "unconditional" here means the schemes have no constraints for the time step only by reason of stability, however, large time step size will definitely induce large errors in practice. This motivates us to develop more accurate schemes, e.g., the second order time marching schemes in this paper.

We must notice that it is very challenging to develop unconditionally energy stable schemes to resolve the stiffness issue induced by the thin interface since the traditional fully-implicit or explicit discretization for the nonlinear term will cause very severe time step constraint (called conditionally energy stable) on the interfacial width [2, 13, 34, 35]. Many efforts had been done (cf. [7, 12, 13, 19, 25, 29, 34, 42, 57] and the references therein) in order to remove this constraint and two commonly used techniques were developed. The first method is the so-called *convex splitting approach* [12, 17, 18, 31, 55, 56], where the convex part of the potential is treated implicitly and the concave part is treated explicitly. The convex splitting approach is unconditionally energy stable, however, it produces the nonlinear scheme, thus the implementation is complicated and the computational cost might be high. The second method is the so-called *stabilized explicit approach* [7, 24, 27, 32–34, 36–39, 42, 43, 45, 46, 53], where the nonlinear term is simply treated explicitly. In order to remove the time step constraint dependence on the interfacial width, a linear stabilizing term has to be added, and the magnitude of that term usually depends on the upper bound of the second order derivative of the double well potential. Such a scheme is efficient and very easy to implement since it is purely linear. Unfortunately, for the second order version, all kinds of linear stabilizers are in vain and the time step constraint still exists (cf. [34]).

Meanwhile, we notice that most of the algorithms developments had been focused on the discretizations for the polynomial type potential, e.g. the double well potential. For the variable mobility Cahn-Hilliard equation with the broadly used logarithmic F-H potential, the efforts about algorithm developments or numerical analysis are even more scarce, although the model had appeared in the original derivation of the Cahn-Hilliard equation [3, 4]. Similar to the double well potential case, the fully-implicit or explicit discretization for the nonlinear term makes the scheme conditionally energy stable as expected, thus it is not efficient in practice. This fact was confirmed in the pioneering work of [8], where the authors developed a temporal first order, nonlinear scheme, where the logarithmic term is discretized in the fully implicit way. Its nonlinear nature can bring up high computational cost as well. Moreover, we emphasize that the two prevalent methods mentioned above are actually not optimal choices for the system (cf. the detailed discussions in Section 3).

Therefore, the main purpose of this paper is to propose some more efficient and accurate time marching schemes for solving the Cahn-Hilliard equation with variable mobility and logarithmic potential. We adopt the *Invariant Energy Quadratization* (IEQ) approach, which is a novel method and has been successfully applied to solve various a variety of phase field type models (cf. [16, 17, 44, 47, 48, 50, 54]). Its essential idea is to transform the bulk potential into a quadratic form (since the nonlinear potential is usually bounded from below) in a set of new variables via the change of variables. Then, for the reformulated model in terms of the new variables that still retains the identical energy dissipation law, all nonlinear terms can be treated semi-explicitly, which then yields a linear system. As a result, we develop a set of efficient schemes which are *accurate* (ready for the second or higher order in time), *easy-to-implement* (linear system), and *unconditionally energy stable* (with a discrete energy dissipation law). We rigorously prove unconditional energy stability for all proposed schemes, including the first order backward Euler, the second order Adam-Bashforth, and the second order Crank-Nicolson schemes. Moreover, we show that the linear operator of all schemes are *symmetric positive definite*, so that one can solve it using the well-developed fast matrix solvers efficiently (CG or other Krylov subspace methods). Through various 2D and 3D numerical simulations, we demonstrate stability and accuracy of the proposed schemes. In addition, although the other typical equation of the phase field approach, i.e., the second order Allen-Cahn equation, is not considered in this paper, all proposed numerical schemes can be easily applied in a similar way.

The rest of the paper is organized as follows. In Section 2, we present the whole system and show the energy law in the continuous level. In Section 3, we develop the numerical schemes and prove their unconditional stabilities. In Section 4, we present various 2D and 3D numerical experiments to validate the accuracy and efficiency of the proposed numerical schemes. Finally, some concluding remarks are presented in Section 5.

## 2. Model Equations

Now we give a brief description for the model equations. We define $\phi$ as the mass concentration of one material component, thus the other component is denoted by $1 - \phi$. The model equation to describe the spinodal decomposition and coarsening phenomena of binary alloys is as follows.

$$\phi_t = \nabla \cdot (M_0(\phi)\nabla\mu), \tag{2.1}$$

$$\mu = -\epsilon^2 \Delta\phi + f(\phi), \tag{2.2}$$

where $M_0(\phi) \geq 0$ is the diffusive mobility function, $f(\phi) = F'(\phi)$ with $F(\phi)$ the Helmholtz free energy, the order parameter $\epsilon$ is the positive constant that measures the interfacial width. For the domain $\Omega$, the boundary condition can be

$$(i) \text{ all variables are periodic; or } (ii) \; \partial_{\mathbf{n}}\phi|_{\partial\Omega} = \nabla\mu \cdot \mathbf{n}|_{\partial\Omega} = 0, \tag{2.3}$$

where $\mathbf{n}$ is the outward normal on the domain boundary $\partial\Omega$.

The thermodynamically reasonable choice is $M_0(\phi) := (\phi - \phi^2)$ (cf. [8]). Such special choice for mobility leads to a number of mathematical difficulties since it is degenerate, i.e., $min(M_0(\phi)) = 0$. Thus, we consider the simpler model following the modified mobility function in [5,8]. For $0 < \sigma \ll 1$, we define the new mobility function $M(\phi)$ to replace $M_0(\phi)$ as follows,

$$M(\phi) = \begin{cases} \frac{1}{4}(1 - (1-\sigma)(2\phi-1)^2), & \text{if } 0 \leq \phi \leq 1, \\ \frac{1}{8}\sigma(1 + e^{\frac{-2(1-\sigma)((2\phi-1)^2-1)}{\sigma}}), & \text{otherwise.} \end{cases} \tag{2.4}$$

It is easy to see that $\exists \, L_1 > L_2 > 0$ such that $L_1 \geq M(x) \geq L_2, \forall x \in \mathbb{R}$.

For the bulk potential, we consider the logarithmic F-H potential, i.e.,

$$\begin{cases} F(\phi) = \phi\ln(\phi) + (1-\phi)\ln(1-\phi) + \theta(\phi - \phi^2), \\ f(\phi) = F'(\phi) = \ln(\dfrac{\phi}{1-\phi}) + \theta(1 - 2\phi), \end{cases} \tag{2.5}$$

where $\theta > 1$.

The model equation (2.1)-(2.2) follows the dissipative energy law. By taking the $L^2$ inner product of (2.1) with $-\mu$, and of (2.2) with $\phi_t$, we have

$$(\phi_t, -\mu) = (M(\phi)\nabla\mu, \nabla\mu), \tag{2.6}$$

$$(\mu, \phi_t) = \frac{d}{dt}\int_\Omega \left(\frac{\epsilon^2}{2}|\nabla\phi|^2 + F(\phi)\right)d\boldsymbol{x}. \tag{2.7}$$

Taking summation of the two equalities, we obtain

$$\frac{d}{dt}E(\phi) = -\|\sqrt{M(\phi)}\nabla\mu\|^2 \leq 0, \tag{2.8}$$

where

$$E(\phi) = \int_\Omega \Big(\frac{\epsilon^2}{2}|\nabla\phi|^2 + F(\phi)\Big)d\boldsymbol{x} \tag{2.9}$$

is the total free energy of the system (2.1)-(2.2).

We remark that the Cahn-Hilliard type dynamical system (2.1)-(2.2) conserves the local mass density. As a matter of fact, by taking the $L^2$ inner product of (2.1) with 1, one can obtain directly

$$\frac{d}{dt}\int_\Omega \phi d\boldsymbol{x} = \int_\Omega \nabla \cdot (M(\phi)\nabla\mu)d\boldsymbol{x} = \int_{\partial\Omega} M(\phi)\partial_\mathbf{n}\mu ds = 0. \tag{2.10}$$

## 3. Numerical Schemes

We now construct several semi-discrete time marching numerical schemes for solving the model system (2.1)-(2.2) and prove their energy stabilities. It will be clear that the energy stabilities of the semi-discrete schemes are also valid in the fully discrete formulation, for instance by finite element or spectral spatial discretizations. To develop the energy stable schemes, the main challenging issues are how to discretize the nonlinear chemical potential term of $f(\phi) = F'(\phi)$ and the mobility function $M(\phi)$.

Following the work in [8], we regularize the logarithmic bulk potential by a $C^2$ piecewise function. More precisely, for any $0 < \sigma \ll 1$, the regularized free energy is

$$\widehat{F}(\phi) = \begin{cases} \phi\ln\phi + \frac{(1-\phi)^2}{2\sigma} + (1-\phi)\ln\sigma - \frac{\sigma}{2} + \theta(\phi-\phi^2), & \text{if} \quad \phi \geq 1-\sigma, \\ \phi\ln\phi + (1-\phi)\ln(1-\phi) + \theta(\phi-\phi^2), & \text{if} \quad \sigma \leq \phi \leq 1-\sigma, \\ (1-\phi)\ln(1-\phi) + \frac{\phi^2}{2\sigma} + \phi\ln\sigma - \frac{\sigma}{2} + \theta(\phi-\phi^2), & \text{if} \quad \phi \leq \sigma. \end{cases} \tag{3.1}$$

For convenience, we consider the problem formulated with the substitute $\widehat{F}(\phi)$, but omit the $\widehat{\phantom{x}}$ in the notation. Now the domain for the regularized functional $F(\phi)$ is now $(-\infty, \infty)$. Hence, we do not need to worry about that any small fluctuation near the domain boundary $(0, 1)$ of the numerical solution can cause the overflow. In [8], the authors proved the error bound between the regularized PDE and the original PDE is controlled by $\sigma$ up to a constant. The regularized energy functional $F(\phi)$ brings up some numerical advantages. First it is convex, thus we can treat it implicitly via the convex splitting approach (cf. [11, 12]). Other alternative method is to treat it explicitly and add some linear stabilizers as the stabilized explicit approach (cf. [34]) if we notice that the second order derivative of $F(\phi)$ is uniformly bounded from above.

Nonetheless, these two techniques are not optimal choices. First, the convex splitting approach brings up nonlinear schemes, thus one need some efficient iterative solvers and the implementations are not easy, which excludes the convex splitting approach. Second, it is shown in [34] that even for the simpler double well potential, the second order version of stabilized explicit method is only conditionally energy stable where the time step is controlled by the interfacial width, which denies the stabilized explicit method.

We aim to develop some more efficient and accurate numerical schemes. In details, we expect that the schemes are linear, unconditionally energy stable, and ready for temporal second order or even higher order. To this end, we adopt the *Invariant Energy Quadratization* (IEQ) approach to design desired numerical schemes, without worrying about whether the continuous/discrete maximum principle holds or a convexity/concavity splitting exists.

3.1. **Transformed system.** We notice that $F(\phi)$ is bounded from below, thus we could rewrite the free energy functional $F(\phi)$ into the following equivalent form

$$F(\phi) = (F(\phi) + B) - B, \tag{3.2}$$

where B is some constant to ensure $F(x) + B > 0, \forall x \in \mathbb{R}$. Therefore, define an auxilliary function $U$ such that

$$U = \sqrt{F(\phi) + B}, \tag{3.3}$$

Thus the total energy of (2.9) turns to a new form

$$E(\phi, U) = \int_\Omega \Big(\frac{\epsilon^2}{2}|\nabla \phi|^2 + U^2\Big) d\boldsymbol{x} - B|\Omega|. \tag{3.4}$$

Then we obtain an equivalent PDE system by taking the time derivative for the new variable $U$:

$$\phi_t = \nabla \cdot (M(\phi)\nabla \mu), \tag{3.5}$$

$$\mu = -\epsilon^2 \Delta \phi + UH, \tag{3.6}$$

$$U_t = \frac{1}{2}H\phi_t, \tag{3.7}$$

where

$$H(\phi) = \frac{f(\phi)}{\sqrt{F(\phi) + B}}, \quad f(\phi) = F'(\phi). \tag{3.8}$$

The boundary conditions for the new system are still (2.3) since the equation (3.7) for the new variable $U$ is simply the ODE with time. The initial conditions read as

$$\phi|_{(t=0)} = \phi_0, \tag{3.9}$$

$$U|_{(t=0)} = \sqrt{F(\phi_0) + B}. \tag{3.10}$$

It is clear that the new transformed system still retains a similar energy dissipative law. Taking the $L^2$ inner product of (3.5) with $-\mu$, of (3.6) with $\phi_t$, of (3.7) with $-U$, and then summing them up, we can obtain the energy dissipation law of the new system (3.5)-(3.6) as

$$\frac{d}{dt}E(\phi, U) = -\|\sqrt{M(\phi)}\nabla \mu\|^2 \leq 0. \tag{3.11}$$

**Remark 3.1.** *We emphasize that the new transformed system (3.5)-(3.7) is exactly equivalent to the original system (2.1)-(2.2), since (3.3) can be easily obtained by integrating (3.7) with respect to the time. For the time-continuous case, the potentials in the new free energy (3.4) are the same as the Lyapunov functional in the original free energy of (2.9), and the new energy law (3.11) for the transformed system is also the same as the energy law (2.8) for the original system as well. We will develop unconditionally energy stable numerical schemes for time stepping of the transformed system (3.5)–(3.6), and the proposed schemes should formally follow the new energy dissipation law (3.11) in the discrete sense, instead of the energy law for the originated system (2.8).*

The time marching numerical schemes are developed to solve the new transformed system (3.5)-(3.6). The proof of the unconditional stability of the schemes follow the similar lines as in the derivation of the energy law (3.11). Let $\delta t > 0$ denote the time step size and set $t^n = n\,\delta t$ for $0 \le n \le N$ with the ending time $T = N\,\delta t$.

## 3.2. First order Scheme. 
We now have the following first order scheme. Having computed $\phi^n, U^n$, we update $\phi^{n+1}, U^{n+1}$ from

$$\frac{\phi^{n+1} - \phi^n}{\delta t} = \nabla \cdot (M(\phi^n)\nabla \mu^{n+1}), \tag{3.12}$$

$$\mu^{n+1} = -\epsilon^2 \Delta \phi^{n+1} + U^{n+1} H^n, \tag{3.13}$$

$$2(U^{n+1} - U^n) = H^n(\phi^{n+1} - \phi^n), \tag{3.14}$$

with the following boundary conditions

$$(i)\ \phi^{n+1}, \mu^{n+1} \text{ are periodic; or } (ii)\ \partial_{\mathbf{n}}\phi^{n+1}|_{\partial\Omega} = \nabla \mu^{n+1} \cdot \mathbf{n}|_{\partial\Omega} = 0. \tag{3.15}$$

and

$$H^n = \frac{f(\phi^n)}{\sqrt{F(\phi^n) + B}}. \tag{3.16}$$

Note that the nonlinear coefficient $H$ of the new variables $U$ are treated explicitly. Moreover, we can rewrite the equations (3.14) as follows:

$$U^{n+1} = \frac{H^n}{2}\phi^{n+1} + (U^n - \frac{H^n}{2}\phi^n), \tag{3.17}$$

Thus (3.12)-(3.13) can be rewritten as the following linear system

$$\frac{\phi^{n+1} - \phi^n}{\delta t} = \nabla \cdot (M(\phi^n)\nabla \mu^{n+1}), \tag{3.18}$$

$$\mu^{n+1} = P_1(\phi^{n+1}) + g_1^n, \tag{3.19}$$

where

$$\begin{cases} P_1(\phi^{n+1}) = -\epsilon^2 \Delta \phi^{n+1} + \frac{1}{2}H^n H^n \phi^{n+1}, \\ g_1^n = U^n - \frac{1}{2}H^n H^n \phi^n. \end{cases} \tag{3.20}$$

Therefore, we can solve $\phi^{n+1}$ and $\mu^{n+1}$ directly from (3.18) and (3.19). Once we obtain $\phi^{n+1}$, the $U^{n+1}$ is automatically given in (3.14). Namely, the new variables $U$ does not involve any extra computational costs. Furthermore, we notice

$$(P_1(\phi), \psi) = \epsilon^2(\nabla \phi, \nabla \psi) + \frac{1}{2}(H^n\phi, H^n\psi), \tag{3.21}$$

if $\psi$ enjoys the same boundary condition as $\phi$ in (3.15). Therefore, the linear operator $P_1(\phi)$ is symmetric (self-adjoint). Moreover, for any $\phi$ with $\int_\Omega \phi d\boldsymbol{x} = 0$, we have

$$(P_1(\phi), \phi) = \frac{\epsilon^2}{2}\|\nabla \phi\|^2 + \frac{1}{2}\|H^n\phi\|^2 \ge 0. \tag{3.22}$$

where " = " is valid if and only if $\phi \equiv 0$.

We first show the well-posedness of the linear system (3.12)-(3.14) (or (3.18)-(3.19)) as follows.

**Theorem 3.1.** *The linear system* (3.18)-(3.19) *admits a unique solution in $H^1(\Omega)$, moreover, the linear operator is symmetric (self-adjoint) and positive definite.*

*Proof.* From (3.18), by taking the $L^2$ inner product with 1, we derive

$$\text{(3.23)} \quad \int_\Omega \phi^{n+1} d\boldsymbol{x} = \int_\Omega \phi^n d\boldsymbol{x} = \cdots = \int_\Omega \phi^0 d\boldsymbol{x}.$$

Let $\alpha_0 = \frac{1}{|\Omega|} \int_\Omega \phi^0 d\boldsymbol{x}$, $\beta_\mu = \frac{1}{|\Omega|} \int_\Omega \mu^{n+1} d\boldsymbol{x}$, and we define

$$\text{(3.24)} \quad \widehat{\phi}^{n+1} = \phi^{n+1} - \alpha_0, \ \widehat{\mu}^{n+1} = \mu^{n+1} - \beta_\mu.$$

Thus, from (3.18)-(3.19), $\widehat{\phi}^{n+1}$ and $\widehat{\mu}^{n+1}$ are the solutions for the following equations,

$$\text{(3.25)} \quad \frac{\psi}{\delta t} - \nabla \cdot (M(\phi^n)\nabla \mu) = f^n,$$

$$\text{(3.26)} \quad \mu + \beta_\mu - P_1(\psi) = \widehat{g}_1^n,$$

with $\widehat{g}_1^n = g_1^n + \frac{1}{2} H^n H^n \alpha_0$, and

$$\text{(3.27)} \quad \int_\Omega \psi d\boldsymbol{x} = 0, \ \int_\Omega \mu d\boldsymbol{x} = 0.$$

We define the inverse quasi-Laplace operator $v = \Delta^{-1} u$ as follows,

$$\text{(3.28)} \quad \begin{cases} \nabla \cdot (M(\phi^n)\nabla v) = u, \\ \int_\Omega v d\boldsymbol{x} = 0, \end{cases}$$

with the periodic boundary condition or Neumann $\partial_{\boldsymbol{n}} v|_{\partial\Omega} = 0$.

Applying $-\Delta^{-1}$ to (3.25) and using (3.26), we obtain

$$\text{(3.29)} \quad -\frac{1}{\delta t}\Delta^{-1}\psi + P_1(\psi) - \beta_\mu = -\Delta^{-1} f^n - \widehat{g}_1^n.$$

Let us express the above linear system (3.29) as $\mathbb{A}\psi = b$,

(i). For any $\psi_1$ and $\psi_2$ in $H^1(\Omega)$ satisfy the boundary conditions (2.3) and $\int_\Omega \psi_1 dx = \int_\Omega \psi_2 dx = 0$, using integration by parts, we derive

$$\text{(3.30)} \quad \begin{aligned} (\mathbb{A}(\psi_1), \psi_2) &= -\frac{1}{\delta t}(\Delta^{-1}\psi_1, \psi_2) + (P_1(\psi_1), \psi_2) \\ &\leq C_1(\|\psi_1\|_{H^{-1}}\|\psi_2\|_{H^{-1}} + \|\nabla\psi_1\|\|\nabla\psi_2\| + \|\psi_1\|\|\psi_2\|) \\ &\leq C_2 \|\psi_1\|_{H^1}\|\psi_2\|_{H^1}. \end{aligned}$$

Therefore, the bilinear form $(\mathbb{A}(\psi_1), \psi_2)$ is bounded $\forall \psi_1, \psi_2 \in H^1(\Omega)$.

(ii). For any $\phi \in H^1(\Omega)$, it is easy to derive that, ,

$$\text{(3.31)} \quad (\mathbb{A}(\phi), \phi) = \frac{1}{\delta t}\|\nabla\Delta^{-1}\phi\|^2 + \frac{\epsilon^2}{2}\|\nabla\phi\|^2 \geq C_3\|\phi\|_{H^1}^2,$$

for $\int_\Omega \phi dx = 0$ from Poincare inequality. Thus the bilinear form $(\mathbb{A}(\phi), \psi)$ is coercive.

Then from the Lax-Milgram theorem, we conclude the linear system (3.29) admits a unique solution in $H^1(\Omega)$.

For any $\psi, \phi$ with $\int_\Omega \psi d\boldsymbol{x} = 0$ and $\int_\Omega \phi d\boldsymbol{x} = 0$, we can easily derive

$$
\begin{aligned}
(\mathbb{A}\psi, \phi) &= (-\frac{1}{\delta t}\Delta^{-1}\psi + P_1(\psi) - \beta_\mu, \phi) \\
&= (\psi, -\frac{1}{\delta t}\Delta^{-1}\phi + P_1(\phi) - \beta_\mu) \\
&= (\psi, \mathbb{A}\phi).
\end{aligned}
\tag{3.32}
$$

Thus $\mathbb{A}$ is self-adjoint. Meanwhile, we derive

$$
\begin{aligned}
(\mathbb{A}\phi, \phi) &= (-\frac{1}{\delta t}\Delta^{-1}\phi + P_1(\phi) - \beta_\mu, \phi) \\
&= \frac{1}{\delta t}\|\nabla\Delta^{-1}\phi\|^2 + (P_1(\phi), \phi) \geq 0,
\end{aligned}
\tag{3.33}
$$

where "=" is valid if only if $\phi = 0$. This concludes the linear operator $\mathbb{A}$ is positive definite. $\square$

The energy stability of the scheme (3.12)-(3.13)-(3.14) (or (3.18)-(3.19)) is presented as follows.

**Theorem 3.2.** *The scheme (3.12)-(3.14) (or (3.18)-(3.19)) is unconditionally energy stable satisfying the following discrete energy dissipation law,*

$$
\frac{1}{\delta t}(E(\phi^{n+1}, U^{n+1}) - E(\phi^n, U^n)) \leq -\|\sqrt{M(\phi^n)}\nabla\mu^{n+1}\|^2 \leq 0,
\tag{3.34}
$$

*where*

$$
E(\phi, U) = \frac{\epsilon^2}{2}\|\nabla\phi\|^2 + \|U\|^2 - B|\Omega|.
\tag{3.35}
$$

*Proof.* By taking the $L^2$ inner product of (3.12) with $-\delta t \mu^{n+1}$, we obtain

$$
-(\phi^{n+1} - \phi^n, \mu^{n+1}) = \delta t \|\sqrt{M(\phi^n)}\nabla\mu^{n+1}\|^2.
\tag{3.36}
$$

By taking the $L^2$ inner product of (3.13) with $\phi^{n+1} - \phi^n$, and applying the following identity

$$
2(a - b, b) = |a|^2 - |b|^2 - |a - b|^2,
\tag{3.37}
$$

we obtain

$$
\begin{aligned}
(\mu^{n+1}, \phi^{n+1} - \phi^n) =& \frac{\epsilon^2}{2}(\|\nabla\phi^{n+1}\|^2 - \|\nabla\phi^n\|^2 + \|\nabla(\phi^{n+1} - \phi^n)\|^2) \\
&+ (U^{n+1}H^n, \phi^{n+1} - \phi^n).
\end{aligned}
\tag{3.38}
$$

By taking the $L^2$ inner product of (3.14) with $-\delta t U^{n+1}$, we obtain

$$
-(\|U^{n+1}\|^2 - \|U^n\|^2 + \|U^{n+1} - U^n\|^2) = -(H^n(\phi^{n+1} - \phi^n), U^{n+1}).
\tag{3.39}
$$

By combining (3.36)-(3.39), we obtain

$$
\begin{aligned}
\frac{\epsilon^2}{2}(\|\nabla\phi^{n+1}\|^2 - \|\nabla\phi^n\|^2 + \|\nabla(\phi^{n+1} - \phi^n)\|^2) + (\|U^{n+1}\|^2 - \|U^n\|^2 + \|U^{n+1} - U^n\|^2) \\
= -\delta t\|\sqrt{M(\phi^n)}\nabla\mu^{n+1}\|^2.
\end{aligned}
\tag{3.40}
$$

Finally, we obtain the desired result (3.34) after dropping some positive terms. $\square$

**Remark 3.2.** *The proposed schemes follow the new energy dissipation law (3.11) formally instead of the energy law for the originated system (2.8). In the time-discrete case, the energy* $E(\phi^{n+1}, U^{n+1})$ *(defined in (3.35)) can be rewritten as a first order approximation to the Lyapunov functionals in* $E(\phi^{n+1})$ *(defined in (2.8)), that can be observed from the following facts, heuristically. From (3.45) and applying the Taylor expansion, we have*

$$(3.41) \qquad U^{n+1} - (\sqrt{F(\phi^{n+1}) + B}) = U^n - (\sqrt{F(\phi^n) + B}) + R_{n+1},$$

*where* $R_{n+1} = O((\phi^{n+1} - \phi^n)^2)$. *Since* $R_k = O(\delta t^2)$ *for* $0 \leq k \leq n+1$ *and* $U^0 = (\sqrt{F(\phi^0) + B})$, *by mathematical induction we can easily get*

$$(3.42) \qquad U^{n+1} = \sqrt{F(\phi^{n+1}) + B} + O(\delta t).$$

**Remark 3.3.** *The essential idea of the IEQ method is to transform the complicated nonlinear potentials into a simple quadratic form in terms set of some new variables via a change of variables.*

*When the nonlinear potential is the double well potential, this method is exactly the same as the so-called Lagrange multiplier method developed in [16]. We remark that the Lagrange multiplier method in [16] only work for the fourth order polynomial potential ($\phi^4$) since its derivative $\phi^3$ can be decomposed into $\lambda(\phi)\phi$ with $\lambda(\phi) = \phi^2$ which can be viewed as a Lagrange multiplier term. However, for other type potentials, the Lagrange multiplier method is not applicable. For example, one cannot separate a factor of $\phi$ from the logarithmic term.*

*On the contrary, the IEQ approach developed here can handle various complex nonlinear terms as long as the corresponding nonlinear potentials are bounded from below. Such an approach does not require the convexity as the convex splitting approach (cf. [12, 17, 18, 31, 55, 56]) or the boundness for the second order derivative as the stabilization approach (cf. [7, 24, 27, 32–34, 36–39, 42, 43, 45, 46, 51–53]). This simple way of quadratization provides some great advantages. First, the complicated nonlinear potential is transferred to a quadratic polynomial form that is much easier to handle. Second, the derivative of the quadratic polynomial is linear, which provides the fundamental support for linearization method. Third, the quadratic formulation in terms of new variables can automatically maintain this property of positivity (or bounded from below) of the nonlinear potentials.*

**Remark 3.4.** *The new variable $U$ is introduced in order to handle the nonlinear bulk potential $F(\phi)$. We notice that the discrete energy still includes the gradient term of $\phi$. Due to the Poincare inequality and mass conservation property for Cahn-Hilliard equation $\int_\Omega \phi^{n+1} d\boldsymbol{x} = \int_\Omega \phi^0 d\boldsymbol{x}$, the $H^1$ stability for the variable $\phi$ is still valid for the proposed scheme, which makes it possible to implement the rigorous error analysis for the proposed schemes. We will implement such analysis in the future work.*

3.3. **Second order Scheme.** We now consider the second-order schemes in time where we discretize the time derivative using the backward differentiation formula (BDF2) and Crank-Nicolson.

3.3.1. *BDF2 Scheme.* Given the initial conditions $\phi^0, U^0$, we first compute $\phi^1, U^1$ by the first order scheme (3.12)-(3.14). Having computed $\phi^n, U^n$ and $\phi^{n-1}, U^{n-1}$, we compute $\phi^{n+1}, U^{n+1}$ as follows,

$$(3.43) \qquad \frac{3\phi^{n+1} - 4\phi^n + \phi^{n-1}}{2\delta t} = \nabla \cdot (M(\phi^*)\nabla \mu^{n+1}),$$

$$(3.44) \qquad \mu^{n+1} = -\epsilon^2 \Delta \phi^{n+1} + U^{n+1} H^*,$$

$$(3.45) \qquad 2(3U^{n+1} - 4U^n + U^{n-1}) = H^*(3\phi^{n+1} - 4\phi^n + \phi^{n-1}),$$

where

$$(3.46) \qquad \begin{cases} \phi^* = 2\phi^n - \phi^{n-1}, \\ H^* = H(\phi^*), \end{cases}$$

and the conditions as follows,

(3.47) $\qquad$ (i) all variables are periodic; or (ii) $\partial_{\mathbf{n}}\phi^{n+1}|_{\partial\Omega} = \nabla\mu^{n+1} \cdot \mathbf{n}|_{\partial\Omega} = 0$,

Similar to the first order schemes, we can rewrite (3.45) as

$$\text{(3.48)} \qquad U^{n+1} = \frac{H^*}{2}\phi^{n+1} + \Big(\frac{4U^n - U^{n-1}}{3} - \frac{H^*}{2}\frac{4\phi^n - \phi^{n-1}}{3}\Big).$$

Thus (3.43)-(3.44) can be rewritten as the following linear system

$$\text{(3.49)} \qquad \frac{3\phi^{n+1} - 4\phi^n + \phi^{n-1}}{2\delta t} = \nabla \cdot (M(\phi^*)\nabla\mu^{n+1}),$$

$$\text{(3.50)} \qquad \mu^{n+1} = P_2(\phi^{n+1}) + g_2^n,$$

where

$$\text{(3.51)} \qquad \begin{cases} P_2(\phi^{n+1}) = -\epsilon^2\Delta\phi^{n+1} + \dfrac{1}{2}H^*H^*\phi^{n+1}, \\ g_2^n = \dfrac{4U^n - U^{n-1}}{3} - \dfrac{1}{2}H^*H^*\dfrac{4\phi^n - \phi^{n-1}}{3}. \end{cases}$$

Therefore, we can solve $\phi^{n+1}$ and $\mu^{n+1}$ directly from (3.49) and (3.50). Once we obtain $\phi^{n+1}$, the $U^{n+1}$ is automatically given in (3.14). Namely, the new variables $U$ does not involve any extra computational costs. Furthermore, we notice

$$\text{(3.52)} \qquad (P_2(\phi), \psi) = \frac{\epsilon^2}{2}(\nabla\phi, \nabla\psi) + \frac{1}{2}(H^*\phi, H^*\psi),$$

if $\psi$ enjoys the same boundary condition as $\phi$ in (3.47). Thus the linear operator $P_2(\phi)$ for the second order scheme is symmetric (self-adjoint). Moreover, for any $\phi$ with $\int_\Omega \phi d\boldsymbol{x} = 0$, we have

$$\text{(3.53)} \qquad (P_2(\phi), \phi) = \frac{\epsilon^2}{2}\|\nabla\phi\|^2 + \frac{1}{2}\|H^*\phi\|^2 \geq 0,$$

where " $=$ " is valid if and only if $\phi = 0$.

We prove the wellposedness and energy stability theorem of scheme (3.43)-(3.47) (or (3.49)-(3.50)) as follows.

**Theorem 3.3.** *The scheme (3.43)-(3.47) (or (3.49)-(3.50)) admits a unique solution in $H^1(\Omega)$, and is unconditionally energy stable satisfying the following discrete energy dissipation law,*

$$\text{(3.54)} \qquad \frac{1}{\delta t}(E_{bdf2}^{n+1} - E_{bdf2}^n) \leq -\|\sqrt{M(\phi^*)}\nabla\mu^{n+1}\|^2 \leq 0,$$

*where*

$$\text{(3.55)} \qquad E_{bdf2}^{n+1} = \frac{\epsilon^2}{2}\Big(\frac{\|\nabla\phi^{n+1}\|^2}{2} + \frac{\|2\nabla\phi^{n+1} - \nabla\phi^n\|^2}{2}\Big) + \Big(\frac{\|U^{n+1}\|^2}{2} + \frac{\|2U^{n+1} - U^n\|^2}{2}\Big) - B|\Omega|.$$

*Moreover, the linear operator is symmetric (self-adjoint) and positive definite.*

*Proof.* The wellposedness and the symmetric positive definite property can be proved similarly as Theorem 3.1.

By taking the $L^2$ inner product of (3.43) with $-2\delta t\mu^{n+1}$, we obtain

$$\text{(3.56)} \qquad -(3\phi^{n+1} - 4\phi^n + \phi^{n-1}, \mu^{n+1}) = 2\delta t\|\sqrt{M(\phi^*)}\nabla\mu^{n+1}\|^2.$$

By taking the $L^2$ inner product of (3.44) with $3\phi^{n+1} - 4\phi^n + \phi^{n-1}$, and applying the following identity

$$(3.57) \qquad (3a - 4b + c, 2a) = a^2 - b^2 + (2a - b)^2 - (2b - c)^2 + (a - 2b + c)^2,$$

we can derive

$$(3\phi^{n+1} - 4\phi^n + \phi^{n-1}, \mu^{n+1})$$
$$(3.58) \qquad = \frac{\epsilon^2}{2}\Big(\|\nabla\phi^{n+1}\|^2 + \|2\nabla\phi^{n+1} - \nabla\phi^n\|^2 - \|\nabla\phi^n\|^2 - \|2\nabla\phi^n - \nabla\phi^{n-1}\|^2\Big)$$
$$+ \frac{\epsilon^2}{2}\|\nabla\phi^{n+1} - 2\nabla\phi^n + \nabla\phi^{n-1}\|^2 + \Big(U^{n+1}H^*, 3\phi^{n+1} - 4\phi^n + \phi^{n-1}\Big).$$

By taking the $L^2$ inner product of (3.45) with $-U^{n+1}$ and applying (3.57), we obtain

$$(3.59) \qquad \begin{aligned}-\Big(\|U^{n+1}\|^2 - \|U^n\|^2 + \|2U^{n+1} - U^n\|^2 - \|2U^n - U^{n-1}\|^2 + \|U^{n+1} - 2U^n + U^{n-1}\|^2\Big)\\ = -\Big(H^*(3\phi^{n+1} - 4\phi^n + \phi^{n-1}), U^{n+1}\Big).\end{aligned}$$

By combinning (3.56), (3.58), and (3.59), we obtain

$$\frac{\epsilon^2}{2}\Big(\|\nabla\phi^{n+1}\|^2 + \|2\nabla\phi^{n+1} - \nabla\phi^n\|^2 - \|\nabla\phi^n\|^2 - \|2\nabla\phi^n - \nabla\phi^{n-1}\|^2\Big)$$
$$(3.60) \quad + \Big(\|U^{n+1}\|^2 + \|2U^{n+1} - U^n\|^2 - \|U^n\|^2 - \|2U^n - U^{n-1}\|^2\Big)$$
$$+ \frac{\epsilon^2}{2}\|\nabla(\phi^{n+1} - 2\phi^n + \phi^{n-1})\|^2 + \|U^{n+1} - 2U^n + U^{n-1}\|^2 = -2\delta t\|\sqrt{M(\phi^*)}\nabla\mu^{n+1}\|^2.$$

Finally, we obtain the desired result (3.54) after dropping some positive terms. □

**Remark 3.5.** *Heuristically, the $\frac{1}{\delta t}(E^{n+1}_{bdf2} - E^n_{bdf2})$ is a second order approximation of $\frac{d}{dt}E(\phi, U)$ at $t = t^{n+1}$. For instance, for any smooth variable $S$ with time, one can write*

$$\Big(\frac{\|S^{n+1}\|^2 + \|2S^{n+1} - S^n\|^2}{2\delta t}\Big) - \Big(\frac{\|S^n\|^2 + \|2S^n - S^{n-1}\|^2}{2\delta t}\Big)$$
$$\cong \Big(\frac{\|S^{n+2}\|^2 - \|S^n\|^2}{2\delta t}\Big) + O(\delta t^2) \cong \frac{d}{dt}\|S(t^{n+1})\|^2 + O(\delta t^2).$$

3.3.2. *Crank-Nicolson Scheme.* Having computed $\phi^n, U^n$, we compute $\phi^{n+1}, U^{n+1}$ as follows.

$$(3.61) \qquad \frac{\phi^{n+1} - \phi^n}{\delta t} = \nabla \cdot (M(\phi^\dagger)\nabla\mu^{n+\frac{1}{2}}),$$

$$(3.62) \qquad \mu^{n+\frac{1}{2}} = -\epsilon^2\Delta\frac{\phi^{n+1} + \phi^n}{2} + \frac{U^{n+1} + U^n}{2}H^\dagger,$$

$$(3.63) \qquad 2(U^{n+1} - U^n) = H^\dagger(\phi^{n+1} - \phi^n),$$

where

$$(3.64) \qquad \begin{cases}\phi^\dagger = \frac{3}{2}\phi^n - \frac{1}{2}\phi^{n-1},\\ H^\dagger = H(\phi^\dagger).\end{cases}$$

and the boundary conditions as follows,

$$(3.65) \qquad (i) \text{ all variables are periodic; or } (ii)\ \partial_\mathbf{n}\phi^{n+1}|_{\partial\Omega} = \nabla\mu^{n+\frac{1}{2}} \cdot \mathbf{n}|_{\partial\Omega} = 0.$$

We note that we can rewrite (3.63) as

$$U^{n+1} = \frac{H^\dagger}{2}\phi^{n+1} + (U^n - \frac{H^\dagger}{2}\phi^n), \tag{3.66}$$

Thus (3.61)-(3.62) can be rewritten as the following linear system

$$\frac{\phi^{n+1} - \phi^n}{\delta t} = \nabla \cdot (M(\phi^\dagger)\nabla \mu^{n+\frac{1}{2}}), \tag{3.67}$$

$$\mu^{n+\frac{1}{2}} = P_3(\phi^{n+1}) + g_3^n, \tag{3.68}$$

where

$$\begin{cases} P_3(\phi^{n+1}) = -\dfrac{\epsilon^2}{2}\Delta\phi^{n+1} + \dfrac{1}{2}H^\dagger H^\dagger \phi^{n+1}, \\ g_3^n = -\dfrac{\epsilon^2}{2}\Delta\phi^n + H^\dagger U^n - \dfrac{1}{4}H^\dagger H^\dagger \phi^n. \end{cases} \tag{3.69}$$

Therefore, we can solve $\phi^{n+1}$ and $\mu^{n+1}$ directly from (3.67) and (3.68). Once we obtain $\phi^{n+1}$, the $U^{n+1}$ is automatically given in (3.63). Namely, the new variables $U$ does not involve any extra computational costs. Furthermore, we notice

$$(P_3(\phi), \psi) = \frac{\epsilon^2}{2}(\nabla\phi, \nabla\psi) + \frac{1}{2}(H^\dagger\phi, H^\dagger\psi), \tag{3.70}$$

if $\psi$ enjoys the same boundary condition as $\phi$ in (3.65). Therefore, the linear operator $P_3(\phi)$ is symmetric (self-adjoint). Moreover, for any $\phi$ with $\int_\Omega \phi d\boldsymbol{x} = 0$, we have

$$(P_3(\phi), \phi) = \frac{\epsilon^2}{2}\|\nabla\phi\|^2 + \frac{1}{2}\|H^\dagger\phi\|^2 \geq 0, \tag{3.71}$$

where " = " is valid if and only if $\phi = 0$.

We prove the wellposedness and energy stability theorem of scheme (3.61)-(3.63) (or (3.67)-(3.68)) as follows.

**Theorem 3.4.** *The scheme (3.61)-(3.63) (or (3.67)-(3.68)) admits a unique solution in $H^1(\Omega)$, and is unconditionally stable satisfying the following discrete energy dissipation law,*

$$\frac{1}{\delta t}(E_{cn2}^{n+1} - E_{cn2}^n) = -\|\sqrt{M(\phi^\dagger)}\nabla\mu^{n+\frac{1}{2}}\|^2 \leq 0, \tag{3.72}$$

*where*

$$E_{cn2}^{n+1} = \frac{\epsilon^2}{2}\|\nabla\phi^{n+1}\|^2 + \|U^{n+1}\|^2 - B|\Omega|. \tag{3.73}$$

*Moreover, the linear operator is symmetric (self-adjoint) and positive definite.*

*Proof.* The wellposedness and the symmetric positive definite property can be proved similarly as Theorem 3.1.

By taking the $L^2$ inner product of (3.61) with $-\delta t \mu^{n+\frac{1}{2}}$, we obtain

$$-(\phi^{n+1} - \phi^n, \mu^{n+\frac{1}{2}}) = \delta t \|\sqrt{M(\phi^\dagger)}\nabla\mu^{n+\frac{1}{2}}\|^2. \tag{3.74}$$

By taking the $L^2$ inner product of (3.62) with $\phi^{n+1} - \phi^n$, we can derive

$$(\phi^{n+1} - \phi^n, \mu^{n+\frac{1}{2}}) = \frac{\epsilon^2}{2}\Big(\|\nabla\phi^{n+1}\|^2 - \|\nabla\phi^n\|^2\Big) + \Big(\frac{U^{n+1} + U^n}{2}H^\dagger, \phi^{n+1} - \phi^n\Big). \tag{3.75}$$

| $\delta t$ | LS1 | Order | LS2-BDF | Order | LS2-CN | Order |
|---|---|---|---|---|---|---|
| $2 \times 10^{-2}$ | 5.87E-2 | - | 3.66E-2 | - | 2.92e-2 | - |
| $1 \times 10^{-2}$ | 4.52E-2 | 0.37 | 1.65E-2 | 1.15 | 1.09e-2 | 1.41 |
| $5 \times 10^{-3}$ | 3.16E-2 | 0.52 | 5.31E-3 | 1.64 | 3.17e-3 | 1.78 |
| $2.5 \times 10^{-3}$ | 1.93E-2 | 0.71 | 1.46E-3 | 1.86 | 8.34e-4 | 1.93 |
| $1.25 \times 10^{-3}$ | 1.06E-2 | 0.86 | 3.78E-4 | 1.95 | 2.12e-4 | 1.98 |
| $6.25 \times 10^{-4}$ | 5.48E-3 | 0.95 | 9.53E-5 | 1.99 | 5.29e-5 | 2.00 |
| $3.125 \times 10^{-4}$ | 2.63E-3 | 1.06 | 2.35E-5 | 2.04 | 1.27e-5 | 2.06 |
| $1.5625 \times 10^{-4}$ | 1.15E-3 | 1.19 | 5.23E-6 | 2.15 | 2.73e-6 | 2.22 |
| $7.8125 \times 10^{-5}$ | 3.84E-4 | 1.58 | 7.68E-7 | 2.76 | 5.97e-7 | 2.19 |

TABLE 4.1. Convergence test for the $L^2$ errors for $\phi$ computed by the first order scheme LS1, the second order scheme LS2-BDF and LS2-CN scheme using different temporal resolutions at $t = 0.5$.

By taking the $L^2$ inner product of (3.63) with $-\frac{U^{n+1}+U^n}{2}$, we obtain

$$(3.76) \qquad -(\|U^{n+1}\|^2 - \|U^n\|^2) = -\Big(H^\dagger(\phi^{n+1} - \phi^n), \frac{U^{n+1} + U^n}{2}\Big).$$

By combining (3.74), (3.75), and (3.76), we obtain

$$(3.77) \qquad \frac{\epsilon^2}{2}\big(\|\nabla\phi^{n+1}\|^2 - \|\nabla\phi^n\|^2\big) + \big(\|U^{n+1}\|^2 - \|U^n\|^2\big) = -\delta t \|\sqrt{M(\phi^\dagger)}\nabla\mu^{n+\frac{1}{2}}\|^2.$$

Finally, we obtain the desired result (3.72). □

**Remark 3.6.** *We emphasize the discrete energy of the Crank-Nicolson scheme strictly follows the PDE energy law (3.11) (i.e., "=" instead of "≤") .*

## 4. Numerical Simulations

We present various 2D and 3D numerical experiments using the schemes constructed in the last section. We set the periodic boundary conditions and use the Fourier-Spectral method to discretize the space. The computational domain is $\Omega = [0, 2\pi]^d, d = 2, 3$. We use $128^d$ Fourier modes in all simulations. The default parameters are

$$(4.1) \qquad \epsilon = 0.05, \sigma = 0.005, B = 3.5, \theta = 2.5.$$

**4.1. Accuracy test.** We first test the convergence rates of the proposed schemes, the first order scheme (3.12)-(3.14) (LS1) and the second order scheme (3.43)-(3.45) (LS2-BDF) and schemes (3.61)-(3.63) (LS2-CN).

We set the initial condition as follows.

$$(4.2) \qquad \phi(t=0) = 0.25\sin(2x)\cos(2y) + 0.48.$$

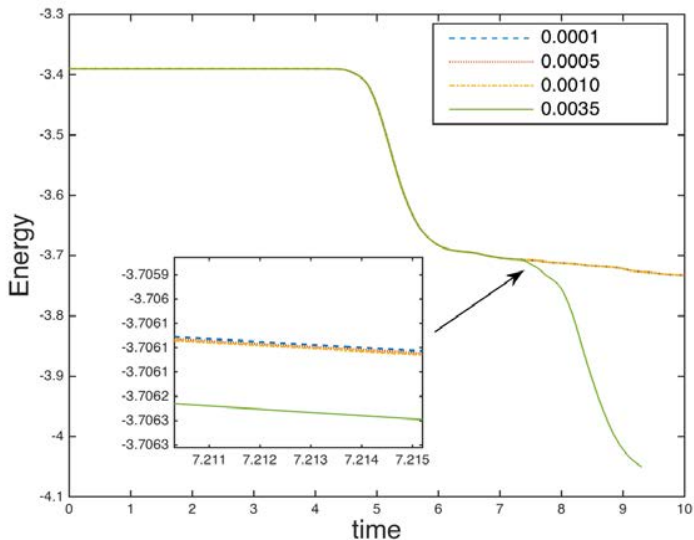

FIGURE 4.1. Time evolution of the free energy functional for four different time steps of $\delta t = 0.0001$, $0.0005$, $0.001$ and $0.0035$ for $\phi_0 = 0.3$ for the 2D case. The energy curves show the decays for all time steps, which confirms that our algorithm is unconditionally stable. The small inset figure shows the small differences in the energy evolution for the three smaller time steps.

Since the exact solutions are not known, we choose the solution obtained by the scheme LS2-CN with the time step size $\delta t = 3.90625 \times 10^{-5}$ as the benchmark solution (approximately the exact solution) for computing errors. We present the $L^2$ errors of the phase variable between the numerical solution and the benchmark solution at $t = 0.5$ with different time step sizes in Table 4.1. We observe that the schemes LS1, LS2-BDF, LS2-CN asymptotically match the first-order and second order accuracy in time, respectively. Moreover, the second order scheme LS2-BDF and LS2-CN gives better and better accuracy along the refinement than the first order scheme LS1 does when using the same time step sizes.

4.2. **Spinodal decomposition in 2D.** In this example, we study the phase separation behavior, the so-called the spinodal decomposition phenomenon using the second order scheme LS2-CN. The process of the phase separation can be studied by considering a homogeneous binary mixture, which is quenched into the unstable part of its miscibility gap. In this case, the spinodal decomposition takes place, which manifests in the spontaneous growth of the concentration fluctuations that leads the system from the homogeneous to the two-phase state. Shortly after the phase separation starts, the domains of the binary components are formed and the interface between the two phases can be specified [1, 9, 57].

The initial condition is taken as the randomly perturbed concentration fields as follows.

(4.3) $$\phi_0 = \bar{\phi} + 0.001 \text{rand}(x, y),$$

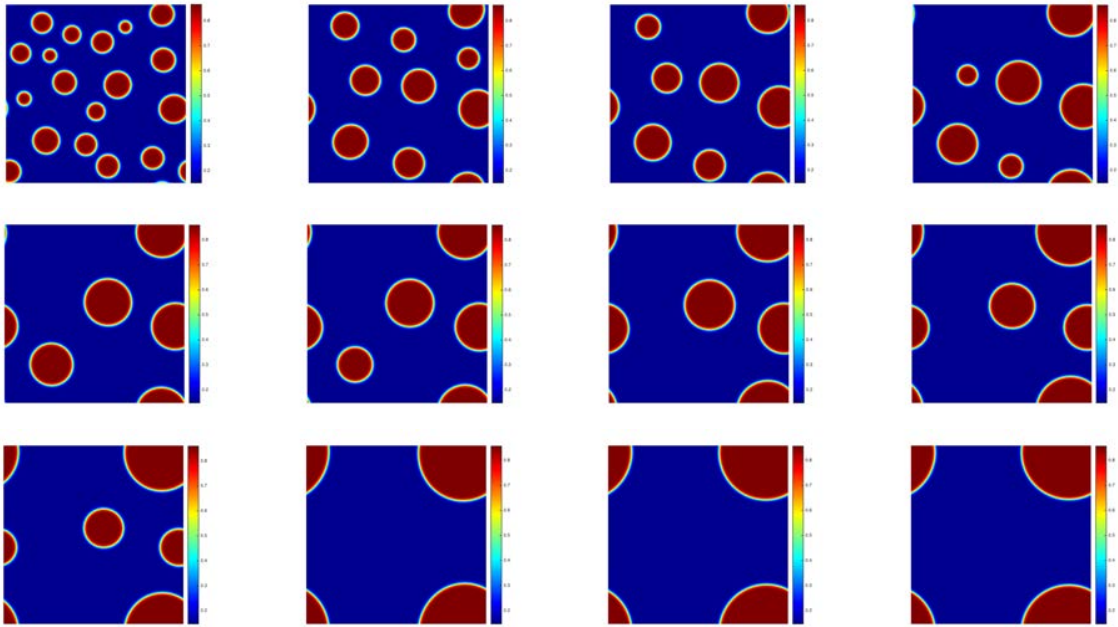

FIGURE 4.2. The 2D dynamical evolution of the local volume fraction phase variable $\phi$ with the initial condition $\bar{\phi}_0 = 0.3$ and time step $\delta t = 0.001$. Snapshots of the numerical approximation are taken at $t = 30, 90, 150, 240, 300, 510, 720, 1000, 1200, 1500, 1710$, and $3000$. The order parameters are of (4.1) and the time step is $\delta t = 0.001$.

where the rand$(x, y)$ is the random number in $[-1, 1]$ and has zero mean. We choose three trial values of $\bar{\phi} = 0.3, 0.5$ and $0.7$ in the following simulations.

We emphasize that any time step size $\delta t$ is allowable for the computations from the stability concern since all developed schemes are unconditionally energy stable. But larger time step will definitely induce large numerical errors. Therefore, we need to discover the rough range of the allowable maximum time step size in order to obtain good accuracy and to consume as low computational cost as possible. Such time step range could be estimated through the energy evolution curve plots, shown in Fig. 4.1, where we compare the time evolution of the free energy for four different time step sizes until $t = 10$ using the second order scheme LS2-CN. We observe that all four energy curves show decays monotonically for all time step sizes, which numerically confirms that our algorithms are unconditionally energy stable. For smaller time steps of $\delta t = 0.0001, 0.0005, 0.001$, the three energy curves coincide very well. But for the larger time step of $\delta t = 0.0035$ ($\sim O(\frac{\epsilon}{2\pi})$), the energy curve deviates viewable away from others. This means the adopted time step size should not be larger than $0.001$, in order to get reasonably good accuracy.

In Fig. 4.2, we perform numerical simulations for initial values $\bar{\phi}_0 = 0.3$ with order parameters (4.1) and time step $\delta t = 0.001$. The red domain, corresponding to the larger values of $\phi$, indicates the concentrated polymer segments [14], and the blue region, corresponding to the smaller values of $\phi$, indicates the macromolecular microspheres (MMs). We observe that the polymer chains are too short to entangle with each other, therefore, they graft on the surface of the MMs. The final equilibrium solution $t = 3000$ presents a very regular shape where the polymer segments accumulate together on the four corners of computed domain.

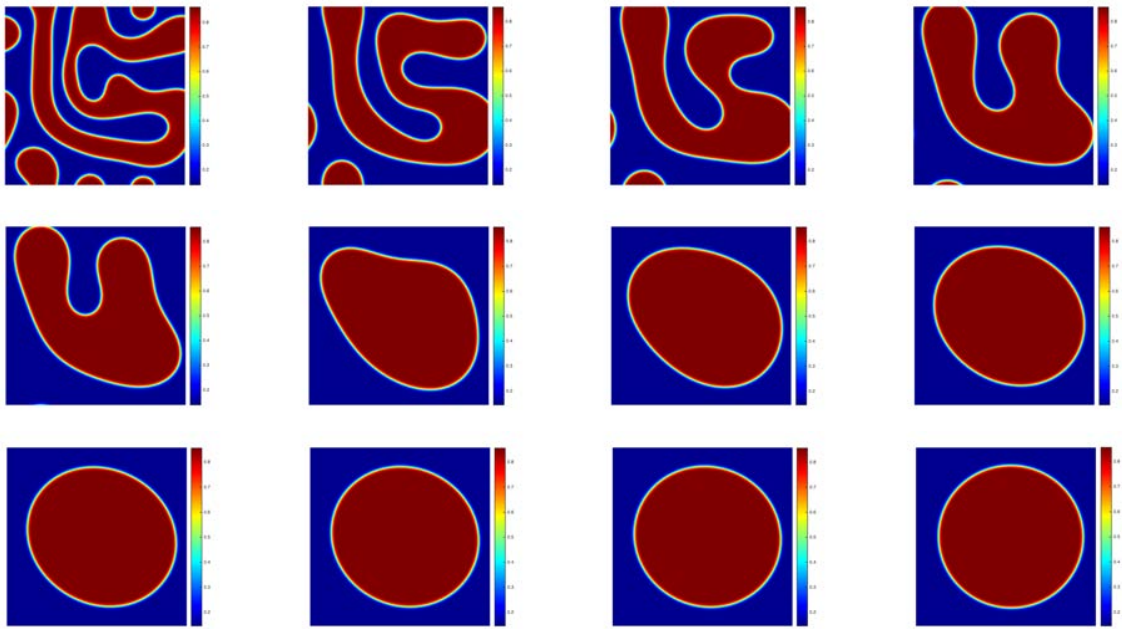

FIGURE 4.3. The 2D dynamical evolution of the local volume fraction phase variable $\phi$ with the initial condition $\bar{\phi}_0 = 0.5$ and time step $\delta t = 0.001$. Snapshots of the numerical approximation are taken at $t = 30, 90, 150, 240, 300, 510, 720, 1000, 1200, 1500, 1710$, and $3000$. The order parameters are of (4.1) and the time step is $\delta t = 0.001$.

We further set the initial values of $\bar{\phi}_0 = 0.5$ in Fig. 4.3 which means the volume fraction of the polymer segments are almost same as the surrounding MMs. The final equilibrium solution is obtained after $t = 3000$, where the polymer segments form a single circular shape. We final set the big value for the volume fraction, i.e., $\bar{\phi}_0 = 0.7$ in Fig. 4.4 which means the volume of the polymer segments are much more than that of the MMs. We observe that the polymer chains entangle firmly and the MMs finally accumulate to circular shape at the four corners. This dynamical behaviors are very similar to the example of $\bar{\phi}_0 = 0.3$ if the colormap is switched.

In Figure 4.5, we present the evolution of the free energy functional for all three initial average of $\phi_0 = 0.3, 0.5$, and $0.7$. The energy curves show the decays with time that confirms that our algorithms are unconditionally stable.

4.3. **Spinodal decomposition in 3D.** We continue to perform the phase separation dynamics using the second order scheme LS2-CN and time step $\delta t = 0.001$, but in 3D space. In order to be consistent with the 2D case, the initial condition reads as follows,

(4.4) $$\phi(t = 0) = \bar{\phi}_0 + 0.001 \mathrm{rand}(x, y, z),$$

where the $\mathrm{rand}(x, y, z)$ is the random number in $[-1, 1]$ with zero mean.

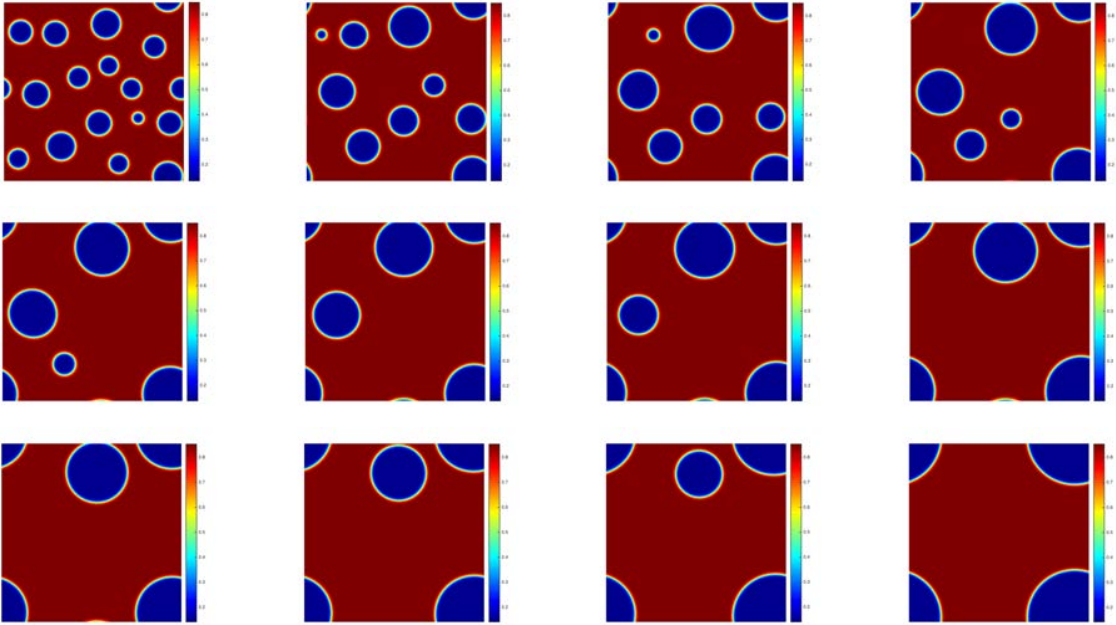

FIGURE 4.4. The 2D dynamical evolution of the local volume fraction phase variable $\phi$ with the initial condition $\bar{\phi}_0 = 0.7$ and time step $\delta t = 0.001$. Snapshots of the numerical approximation are taken at $t = 30, 90, 150, 240, 300, 510, 720, 1000, 1200, 1500, 1710$, and $3000$. The order parameters are of (4.1) and the time step is $\delta t = 0.001$.

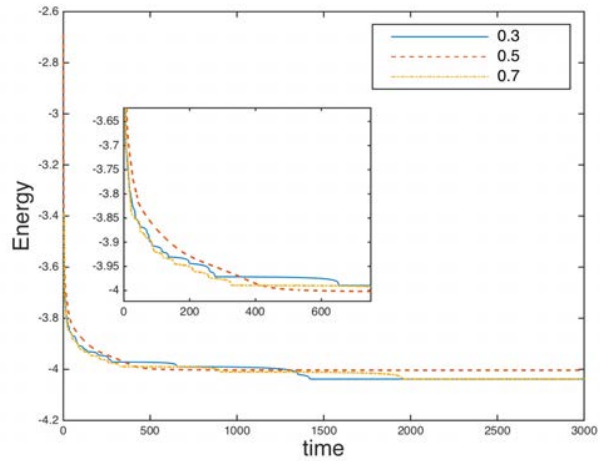

FIGURE 4.5. Time evolution of the free energy functional using the algorithm LS2-CN using $\delta t = 0.001$ for the three initial values of $\bar{\phi}_0 = 0.3, 0.5$ and $0.7$ for the 2D case.

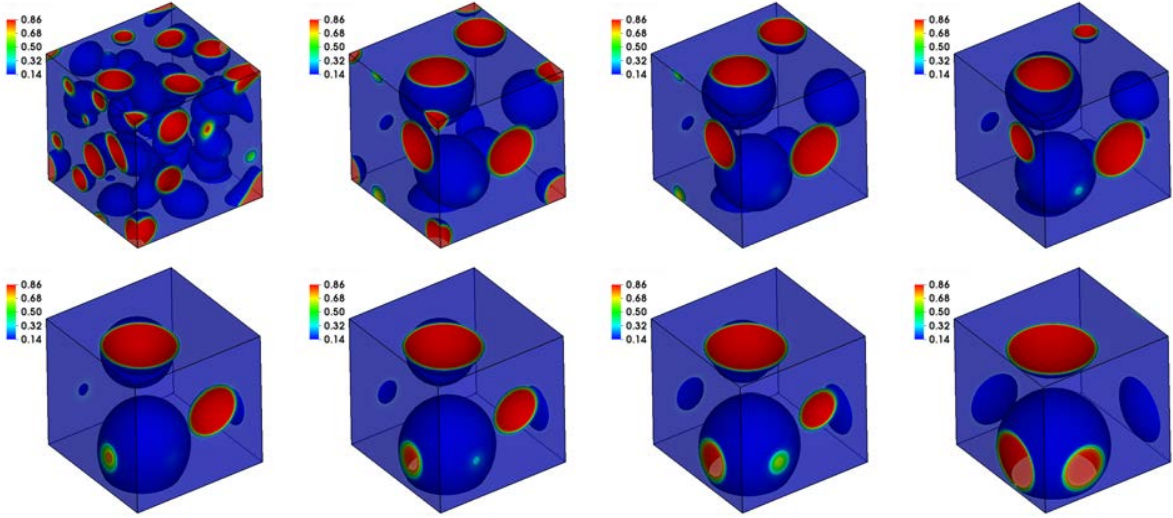

FIGURE 4.6. The 3D dynamical evolution of the local volume fraction phase variable $\phi$ with the initial condition $\bar{\phi}_0 = 0.3$ and time step $\delta t = 0.001$. Snapshots of the numerical approximation are taken at $t = 50, 250, 350, 500, 900, 1000, 1100$ and $2500$. The order parameters are of (4.1) and the time step is $\delta t = 0.001$.

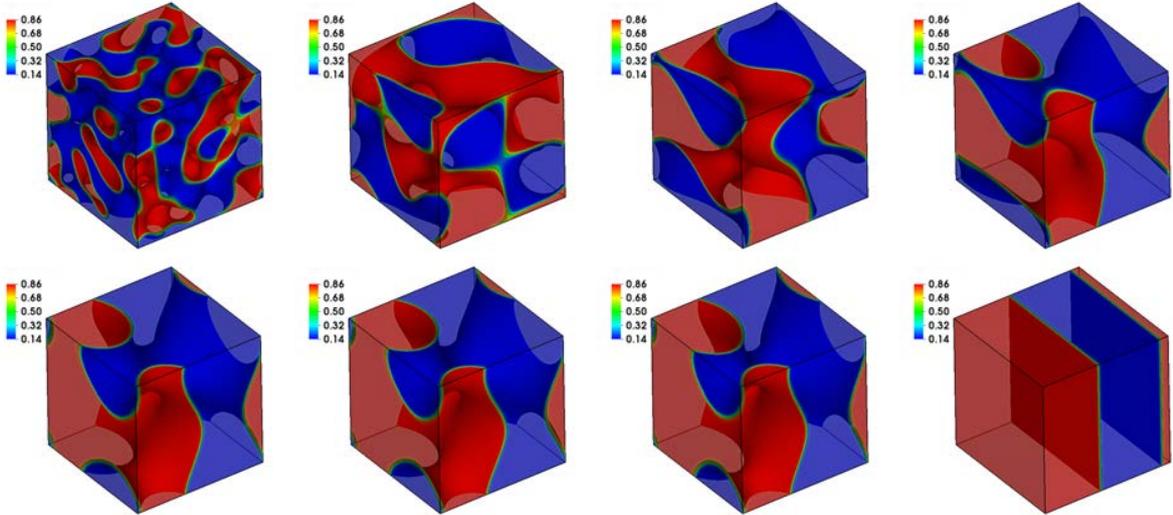

FIGURE 4.7. The 3D dynamical evolution of the local volume fraction phase variable $\phi$ with the initial condition $\bar{\phi}_0 = 0.5$ and time step $\delta t = 0.001$. Snapshots of the numerical approximation are taken at $t = 50, 500, 1000, 1500, 2500, 3500, 5000$ and $7500$. The order parameters are of (4.1) and the time step is $\delta t = 0.001$.

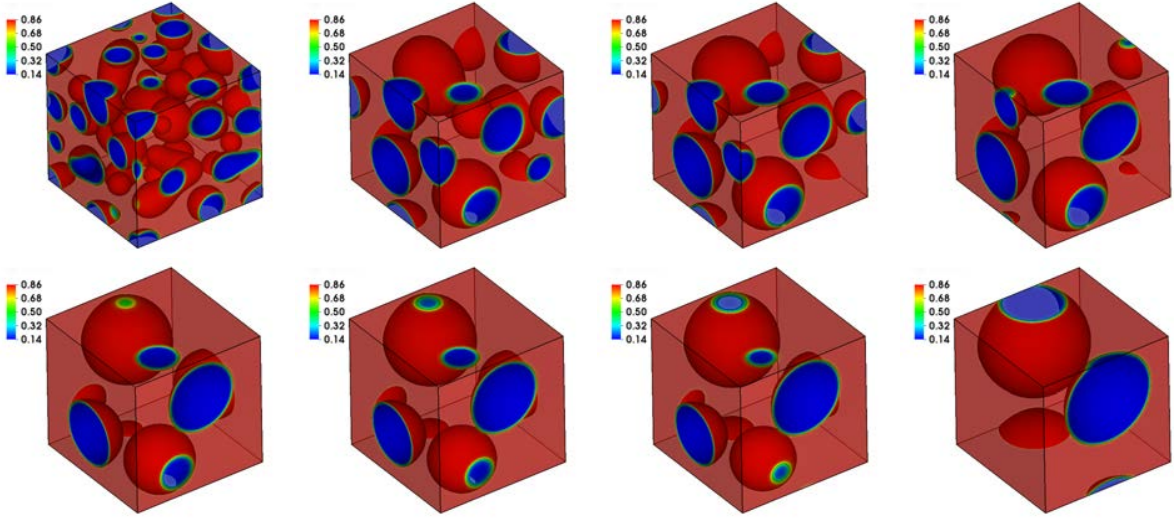

FIGURE 4.8. The 3D dynamical evolution of the local volume fraction phase variable $\phi$ with the initial condition $\bar{\phi}_0 = 0.7$ and time step $\delta t = 0.001$. Snapshots of the numerical approximation are taken at $t = 50, 250, 350, 500, 900, 1000, 1100$ and $2500$. The order parameters are of (4.1) and the time step is $\delta t = 0.001$.

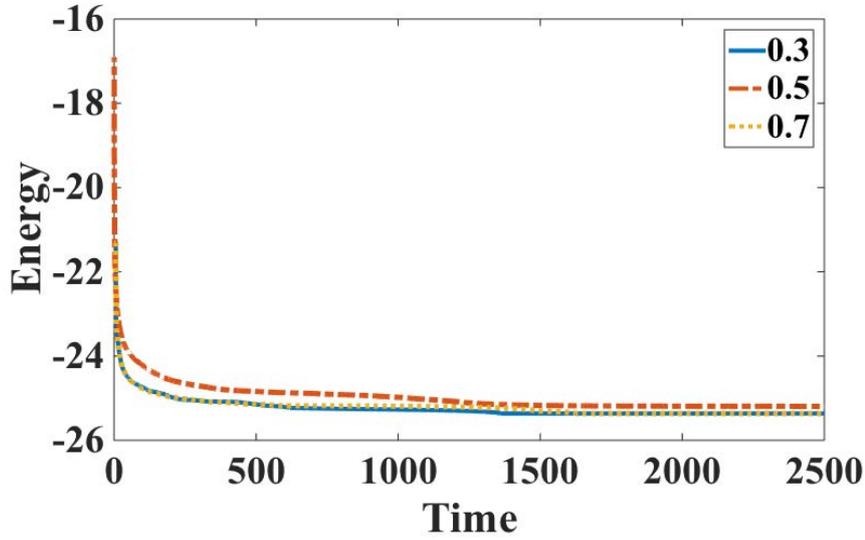

FIGURE 4.9. Time evolution of the free energy functional using the algorithm LS2-CN using $\delta t = 0.001$ for the three initial values of $\bar{\phi}_0 = 0.3, 0.5$ and $0.7$ in 3D case.

Fig. 4.6 shows the dynamical behaviors of the phase separation for the initial value $\bar{\phi}_0 = 0.3$. The dynamical evolutions are consistent with the 2D example Fig. 4.2. We observe that the polymer segments finally accumulate together to the spherical shape, as shown in 2D case. For initial value $\bar{\phi}_0 = 0.5$, in Fig. 4.7, we observe that the polymer segments initially partially accumulate as shown in 2D case. However, the final steady state forms the uniform two layer structure. In Fig. 4.8, we set the initial value $\bar{\phi}_0 = 0.7$. We observe that the polymers chains firmly entangled with each other and the MMs finally form the spherical shape, as Fig. 4.4 of the 2D case. In Fig. 4.9, for all three initial average, we present the evolution of the free energy functional, in which the energy curves show the decays with time.

## 5. Concluding Remarks

In this paper, we develop three efficient schemes to solve the variable mobility Cahn-Hilliard system with the logarithmic Flory-Huggins potential, using a novel, so called IEQ approach. Such approach is effective and efficient, and particularly suitable to discretize the complicated nonlinear potential that is bounded from below. Compared to the prevalent nonlinear schemes based on the convex splitting approaches or other nonlinear schemes, the IEQ approach can easily conquer the inconvenience from nonlinearities by linearizing the nonlinear terms in the new way. The developed schemes (i) are *accurate* (ready for second or higher order in time); (ii) are *stable* (unconditional energy dissipation law holds); and (iii) are *easy to implement* (only need to solve linear equations at each time step). Furthermore, the induced linear system is symmetric positive definite, thus one can apply any Krylov subspace methods with mass lumping as pre-conditioners for solving such system efficiently. We emphasize that, to the best of the authors' knowledge, the schemes to solve the case of logarithmic potential are the first such linear and accurate schemes with provable energy stabilities.

Finally, the method is general enough to be extended to develop linear schemes for a large class of gradient flow problems with complex nonlinearities in the free energy density. Although we consider only time discrete schemes in this study, the results can be carried over to any consistent finite-dimensional Galerkin approximations since the proofs are all based on a variational formulation with all test functions in the same space as the space of the trial functions.

**Acknowledgments.** X. Yang is partially supported by NSF DMS-1200487 and NSF DMS-1418898.